\documentclass[final,1p,times]{elsarticle}
\usepackage{amssymb}
\usepackage{amsfonts}
\usepackage{amsmath}
\usepackage{geometry}
\usepackage{graphicx}
\usepackage{epsfig}
\usepackage{epstopdf}

\usepackage{bm}
\newmuskip\pFqmuskip

\newcommand*\pFq[6][8]{%
  \begingroup 
  \pFqmuskip=#1mu\relax
  \mathchardef\normalcomma=\mathcode`,
  \mathcode`\,=\string"8000
  \begingroup\lccode`\~=`\,
  \lowercase{\endgroup\let~}\pFqcomma
  {}_{#2}F_{#3}{\left[\genfrac..{0pt}{}{#4}{#5};#6\right]}%
  \endgroup
}
\newcommand{\pFqcomma}{{\normalcomma}\mskip\pFqmuskip}
\numberwithin{equation}{section}
\newtheorem{thm}{Theorem}

\newdefinition{rmk}{Remark}

\newproof{pf}{Proof} 

\newproof{pot}{Proof of Theorem  
}
\begin{document}
\begin{frontmatter}
\title{Impossibility of convergence of a confluent Heun function on the boundary of the disc of convergence}
\author{Yoon-Seok Choun}
\ead{ychoun@gradcenter.cuny.edu; ychoun@gmail.com}
\address{Department of Physics, Hanyang University, Seoul, 133-791, South Korea}
\begin{abstract}
 
 The confluent Heun equation is one of 4 confluent forms of Heun's differential equation in which is the Fuchsian equation of second order with four regular singularities. A confluent Heun function is applicable to diverse areas such as theory of rotating/non-rotating black hole, the gauge theories on thick brane words, Schr$\ddot{\mbox{o}}$dinger equation for hydrogen molecule ion in Stark effect and etc.
 The confluent Heun function consists of the three term recurrence relation in its power series, and we show that the function is divergent on the boundary of the disc of convergence. 
\end{abstract}

\end{frontmatter}
\section{Introduction}\label{sec.1} 
In 1889, Karl L. W. Heun suggested the Fuchsian equation of second order with four regular singularities such as 
\begin{equation}
\frac{d^2{y}}{d{x}^2} + \left(\frac{\gamma }{x} +\frac{\delta }{x-1} + \frac{\varepsilon }{x-a}\right) \frac{d{y}}{d{x}} +  \frac{\alpha \beta x-q}{x(x-1)(x-a)} y = 0 \nonumber
\end{equation}
with the condition $\varepsilon = \alpha +\beta -\gamma -\delta +1$ and $a \ne 0 $ \cite{Heun1889,Ronv1995}. Currently, its equation is called as the Heun equation and its equation has the four kinds of confluent forms such as  confluent Heun, doubly confluent Heun,  biconfluent Heun and triconfluent Heun equations. 
We can derive these confluent forms from the Heun equation by combining two or more regular singularities to take form an irregular singularity: Its process is similar to deriving of confluent hypergeometric equation from the hypergeometric equation.

The non-symmetrical canonical form of the Confluent Heun Equation is taken as \cite{Deca1978,Decar1978,Ronv1995}  
\begin{equation}
\frac{d^2{y}}{d{x}^2} + \left(\beta  +\frac{\gamma }{x} + \frac{\delta }{x-1}\right) \frac{d{y}}{d{x}} +  \frac{\alpha \beta x-q}{x(x-1)} y = 0 \label{eq:1}
\end{equation}
(\ref{eq:1}) has three singular points: two regular singular points which are 0 and 1 with exponents $\{0, 1-\gamma\}$ and $\{0, 1-\delta \}$, and one irregular singular point which is $\infty$ with an exponent $\alpha$. If $\beta =0$, (\ref{eq:1}) is the hypergeometric equation well known to us. In this paper, we assume $\beta \ne 0$.

Assume that $y(x)$ has a series expansion of the form
\begin{equation}
y(x)= \sum_{n=0}^{\infty } d_n x^{n+\lambda } \label{eq:2}
\end{equation}
where $\lambda $ is an indicial root. Plug (\ref{eq:2})  into (\ref{eq:1}):
\begin{equation}
d_{n+1}=A_n \;d_n +B_n \;d_{n-1} \hspace{1cm};n\geq 1 \label{eq:3}
\end{equation}
with  
\begin{equation}
\begin{cases} 
A_n  = \dfrac{n^2 + (2\lambda -\beta +\gamma +\delta -1)n +\lambda (\lambda -\beta +\gamma +\delta -1)-q}{n^2 +(2\lambda +1+\gamma  )n +(\lambda +1)(\lambda +\gamma )} \\
B_n   = \beta \dfrac{  n +\lambda +\alpha -1}{n^2 +(2\lambda +1+\gamma  )n +(\lambda +1)(\lambda +\gamma )} \\
d_1= A_0 \;d_0  
\end{cases}\label{eq:4}
\end{equation}
We have two indicial roots which are $\lambda = 0$ and $ 1-\gamma $. And a confluent Heun function of (\ref{eq:1}) around $x=0$ is absolutely convergent where $\left| x\right| <1 $. 

One of well-known examples of the 2-term recurrence relation in a power series is a hypergeometric function. And its domain is given by 
\begin{equation}
\mathcal{D} :=  \left\{  x \in \mathbb{C} \big|  \lim_{n\rightarrow \infty}\left| \frac{(n+a)(n+b)}{(n+c)(n+1)} x\right|  = \left|  x \right| <1 \right\} \nonumber
\end{equation}
where $a,b,c \in \mathbb{C}$. 

In 1812 Carl Friedrich Gauss published a definitive test that determines convergence for a hypergeometric series at $|x|=1$, and its series is absolute convergent as $\mathbb{R}(c) > \mathbb{R}(a+b)$ \cite{Birk,Gauss}.  
But it has not been known which value coefficients for the 3-term recursive relation in a power series on the boundary of the disc of convergence make the series as absolute convergent. To put it another way, what condition of coefficients makes a series as absolutely convergent at $ \left| \lim_{n\rightarrow \infty} A_n x \right| =|x|=1 $ where $ \left| \lim_{n\rightarrow \infty} B_n \right| =0 $ for a confluent Heun function?  
In this paper, we show why  the function is not convergent on the boundary of the disc of convergence.

\section{Result}\label{sec.2}
From a methodological point of view of the demonstration, the nature of a confluent Heun function is studied by utilizing Gauss' test in a similar way as Gauss studied a hypergeometric series on the boundary of the disc of convergence.   
\begin{thm} 
Let $y(x)= x^{\lambda }\left( d_0 +d_1 x+ d_2 x^2+ d_3 x^3+ \cdots \right)$, $d_i \ne 0$, be a power series (including a confluent Heun function) which consists of the 3-term recurrence relation. We assume that the polynomials in the numerator and denominator of $A_n$ have the same degree, and the degree of polynomials in the numerator of $B_n$ is  smaller than the degree of polynomials in the denominator of it:
\begin{equation}
d_{n+1}=A_n \;d_n +B_n \;d_{n-1} \hspace{1cm};n\geq 1 \label{eq:77}
\end{equation} 
\begin{equation}
\begin{cases} 
A_n =  \dfrac{C_t n^t + C_{t-1}n^{t-1} + \cdots + C_0}{c_t n^t + c_{t-1}n^{t-1} + \cdots + c_0}  = A  \dfrac{ n^t + \Omega _{t-1}n^{t-1} + \cdots + \Omega _0}{ n^t + \omega _{t-1}n^{t-1} + \cdots + \omega _0} \\
B_n = \dfrac{ G_{t-1}n^{t-1} + G_{t-2}n^{t-2} + \cdots + G_0}{g_t n^t + g_{t-1}n^{t-1} + \cdots + g_0}  = B \dfrac{ n^{t-1} + \Theta _{t-2}n^{t-2} + \cdots + \Theta _0}{ n^t + \theta _{t-1}n^{t-1} + \cdots + \theta _0} \\
d_1= \dfrac{C_0}{c_0}  \;d_0 = A \dfrac{\Omega _0}{\omega _0}  \;d_0   
\end{cases}\nonumber
\end{equation}
where neither $C_t$ nor $c_t$ is zero and neither $G_{t-1}$ nor $g_t$ is zero with $\Omega _j= C_j/C_t$, $\omega _j = c_j/c_t $, $\Theta _j = G_j/G_{t-1} $ and $\theta _j= g_j/g_t$. Also, we denote $A=C_t/c_t$ and $B= G_{t-1}/g_t$.
 
The domain of absolute convergence of $y(x)$ is given by \cite{Choun2013} 
\begin{equation}
\mathcal{D} :=  \left\{  x \in \mathbb{C} \Bigg|  \left| \lim_{n\rightarrow \infty} A_n x \right| + \left| \lim_{n\rightarrow \infty} B_n x^2 \right| =  \left| A x \right| <1 \right\}
 \nonumber
\end{equation}
If $ \left| A x \right| =1$, then the series cannot converges. \label{Thm.1} 
\end{thm} 
\begin{pf}
\textbf{1. The case of} $\mathbf{\Omega _{t-1}< \omega _{t-1}}$ 

We define
\begin{equation} 
\left| \overline{A}_n \right| =  \left|   \dfrac{ n^t + \Omega _{t-1}n^{t-1} + \cdots + \Omega _0}{ n^t + \omega _{t-1}n^{t-1} + \cdots + \omega _0}  \right|  \;\;\mbox{and}\;\; \left| \overline{B}_n \right| =  \left|  \dfrac{ n^{t-1} + \Theta _{t-2}n^{t-2} + \cdots + \Theta _0}{ n^t + \theta _{t-1}n^{t-1} + \cdots + \theta _0}  \right|    
 \nonumber
\end{equation}  
We say $W(n)= n^t + \Omega _{t-1}n^{t-1} + \cdots + \Omega _0$ and $w(n)= n^t + \omega _{t-1}n^{t-1} + \cdots + \omega _0$.
  
For  a large $n$, determining the sign of the polynomial is the leading coefficient.
 After passing the rightmost root of a polynomial, a polynomial takes only positive values if  the coefficient of the coefficient of the highest power of $n$ is positive \cite{Bress}. Once $n$ is larger than the largest root of $W(n)$ or $w(n)$, then 
\begin{equation} 
\left| \overline{A}_n \right| =   \frac{W(n)}{w(n)}
 \nonumber
\end{equation}  
If $\Omega _{t-1}< \omega _{t-1}$, then we can find a positive integer $h$ such that $\omega _{t-1} - \Omega _{t-1}- h <0 $. Then, for a large value of $n$, we know 
\begin{equation}
\frac{n}{n-h} \left| \overline{A}_n \right| = \dfrac{ n^{t+1} + \Omega _{t-1}n^{t} + \cdots }{ n^{t+1} + \left( \omega _{t-1}- h \right) n^{t} + \cdots } >1
\nonumber
\end{equation}
because of $ \omega _{t-1}- h < \Omega _{t-1}$. So we have
\begin{equation}
 \left| \overline{A}_n \right| > 1-\frac{h}{n}
\nonumber
\end{equation}
Given  error bound $0< \epsilon \ll 1$ and a positive integer $N$,
$n\geq N$ implies that 
\begin{equation}
 \left| \overline{A}_n \right| > 1-\frac{h}{n} > 1-\epsilon 
\label{eq:5}
\end{equation}
For a large value of $n$,
\begin{equation}
\left| \overline{B}_n \right| \sim \frac{1}{n}
\nonumber
\end{equation}
Then, we arrive at
\begin{equation}
 \left| \overline{B}_n \right| > \frac{1}{n+h_0}
\label{eq:6}
\end{equation}
for a large $n$ and some positive integer $h_0$.

We express $\overline{i}= N+i$ and $\hat{i}= N+1+i$ where $i\in \{ 0,1,2,\cdots \}$.
For $n=N, N+1, N+2, \cdots$ in succession, take the modulus of the general term of $d_{n+1}$ in  (\ref{eq:77})
\begin{equation} 
 |d_{N+j}|= |\overline{c}_{j}| |d_N| + |\hat{c}_{j-1}| |B_{\overline{0}}|  |d_{N-1}|   
  \label{eq:7}
\end{equation}
where $j\in \{1,2,3,\cdots \}$, and  we define $|\overline{c}_{0}|= |\hat{c}_{0}|=1$.
 
In (\ref{eq:7}) is $|\overline{c}_{n}|$ is the sequence of the 3-term recurrence relation such as 
\begin{equation} 
 |\overline{c}_{n+1}|= \left|A_{\overline{n}} \right|  |\overline{c}_{n}| + \left|B_{\overline{n}} \right|  |\overline{c}_{n-1}|,\;\; n\geq 1   
  \label{eq:10}
\end{equation}
where $ |\overline{c}_{1}|= \left|A_{\overline{0}} \right|$.

Similarly, $|\hat{c}_{n}|$ is the sequence of the 3-term recurrence relation such as
\begin{equation} 
 |\hat{c}_{n+1}|= \left|A_{\hat{n}} \right|  |\hat{c}_{n}| + \left|B_{\hat{n}} \right|  |\hat{c}_{n-1}|,\;\; n\geq 1   
  \label{eq:11}
\end{equation}
where $ |\hat{c}_{1}|= \left|A_{\hat{0}} \right|$.

According to (\ref{eq:7}), then the series of absolute values, $1+|d_1||x|+|d_2||x|^2+ |d_3||x|^3 +\cdots$, is dominated by the convergent series
\begin{multline}
\sum_{n=0}^{N-1}|d_n||x|^n + |d_N||x|^N + \left( \overline{c}_{1}| |d_N| +   |B_{\overline{0}}|  |d_{N-1}|  \right) |x|^{N+1}
+ \left( \overline{c}_{2}| |d_N| + |\hat{c}_{1}| |B_{\overline{0}}|  |d_{N-1}| \right) |x|^{N+2} \\
+\left( \overline{c}_{3}| |d_N| + |\hat{c}_{2}| |B_{\overline{0}}|  |d_{N-1}|  \right) |x|^{N+3}
+\left( \overline{c}_{4}| |d_N| + |\hat{c}_{3}| |B_{\overline{0}}|  |d_{N-1}|  \right) |x|^{N+4}+\cdots \\
= \sum_{n=0}^{N-1}|d_n||x|^n +   |d_N|  |x|^N \sum_{i=0}^{\infty }\left| \overline{c}_{i} \right||x|^{i} + |B_{\overline{0}}|  |d_{N-1}|  |x|^{N+1} \sum_{i=0}^{\infty }\left| \hat{c}_{i}\right| |x|^{i}
\label{eq:12}
\end{multline}
By rearranging coefficients $\left| A_{\overline{n}}\right|$ and $\left|B_{\overline{n}}\right|$ in each sequence $\left|\overline{c}_{n}\right|$ in (\ref{eq:10}), 
\begin{equation}
\sum_{i=0}^{\infty }\left|\overline{c}_{i}\right| |x|^{i} = \left| \overline{y}_0(z)\right| + \left|\overline{y}_1(z)\right| \eta +\sum_{\tau =2}^{\infty } \left|\overline{y}_{\tau }(z)\right| \eta ^{\tau }  \label{eq:13}
\end{equation}
where
\begin{align}
\left|\overline{y}_0(z)\right| &= \sum_{i_0=0}^{\infty } \prod _{i_1=0}^{i_0-1} \left| \overline{A}_{i_1+N} \right| z^{i_0  }\nonumber\\
\left|\overline{y}_1(z) \right| &=  \sum_{i_0=0}^{\infty }  \left| \overline{B}_{i_0+1+N}\right| \prod _{i_1=0}^{i_0-1} \left| \overline{A}_{i_1+N} \right| \sum_{i_2=i_0}^{\infty }  \prod _{i_3=i_0}^{i_2-1} \left| \overline{A}_{i_3+2+N}\right|  z^{i_2 }\nonumber\\
\left|\overline{y}_{\tau }(z)\right|  &=  \sum_{i_0=0}^{\infty } \left| \overline{B}_{i_0+1+N}\right| \prod _{i_1=0}^{i_0-1} \left|\overline{A}_{i_1+N} \right|
\prod _{k=1}^{\tau -1} \Bigg( \sum_{i_{2k}= i_{2(k-1)}}^{\infty } \left|\overline{B}_{i_{2k}+2k+1+N}\right| \prod _{i_{2k+1}=i_{2(k-1)}}^{i_{2k}-1} \left|\overline{A}_{i_{2k+1}+2k+N}\right|\Bigg)\nonumber\\
  & \hspace{.3cm}\times \sum_{i_{2\tau } = i_{2(\tau -1)}}^{\infty }   \prod _{i_{2\tau +1}=i_{2(\tau -1)}}^{i_{2\tau }-1} \left|\overline{A}_{i_{2\tau +1}+2\tau +N}\right| z^{i_{2\tau } }    \label{eq:14}
\end{align} 
and
\begin{equation}
\begin{cases}
\eta = |B| |x|^2 \cr
z= |A| |x|    
\end{cases}\nonumber 
\end{equation}
 The sequence $\left|\overline{c}_{n}\right|$ combines into combinations of $\left|A_{\overline{n}}\right|$ and $\left|B_{\overline{n}}\right|$ terms in (\ref{eq:10}): (\ref{eq:13}) is done by letting $\left|B_{\overline{n}}\right|$ in the sequence $\left|\overline{c}_{n}\right|$ is the leading term in a series $\sum_{i=0}^{\infty }\left|\overline{c}_{i}\right| |x|^{i}$; we observe the term of sequence $\left|\overline{c}_{n}\right|$ which includes zero term of $\left|B_{\overline{n}}\right|'s$ for a sub-power series $\left|\overline{y}_0(z)\right|$, one term of $\left|B_{\overline{n}}\right|'s$ for the sub-power series $\left|\overline{y}_1(z)\right|$, two terms of $\left|B_{\overline{n}}\right|'s$ for a $\left|\overline{y}_2(z)\right|$, three terms of $\left|B_{\overline{n}}\right|'s$ for a $\left|\overline{y}_3(z)\right|$, etc.

Similarly, by rearranging coefficients $\left| A_{\hat{n}}\right|$ and $\left|B_{\hat{n}}\right|$ in each sequence $\left|\hat{c}_{n}\right|$ in (\ref{eq:11}), 
\begin{equation}
\sum_{i=0}^{\infty }\left|\hat{c}_{i}\right| |x|^{i} = \left| \hat{y}_0(z)\right| + \left|\hat{y}_1(z)\right| \eta +\sum_{\tau =2}^{\infty } \left|\hat{y}_{\tau }(z)\right| \eta ^{\tau }  \label{eq:15}
\end{equation}
here, $\left| \hat{y}_k(z)\right|$ where $k\in \{0,1,2,\cdots\}$ is given by replacing an index $N$ with $N+1$ in (\ref{eq:14}).

The asymptotic of the Pochhammer symbol $(a)_n$ as $n\rightarrow \infty $ is given by
\begin{equation}
(a)_n \sim \frac{\sqrt{2 \pi}}{\Gamma (a)} e^{-n} n^{n+a-1/2} \left( 1+ \mathcal{O}(1/n)\right)
\nonumber
\end{equation}
and
\begin{equation}
\frac{(a)_n}{(b)_n} \sim \frac{\Gamma (b)}{\Gamma (a)}  n^{a-b} \left( 1+ \mathcal{O}(1/n)\right), \hspace{1cm} n\rightarrow \infty 
\nonumber
\end{equation} 
then, we can say
\begin{eqnarray}
\frac{\left( N+2(r+1)-h +i_{2r}\right)_{i_{2(r+1)}}}{\left(  N+2(r+1) +i_{2r}\right)_{i_{2(r+1)}}} 
&>& \frac{\Gamma \left(  N+2(r+1) +i_{2r} \right)}{2\Gamma \left(  N+2(r+1)-h +i_{2r}\right)}  \frac{1}{ i_{2(r+1)}^{h}} 
\label{eq:16}
\end{eqnarray} 
where $i_{2(r+1)} \geq m$, here $N$ and $h$ are some positive integers, $r\in \{1,2,3,\cdots\}$ and  $N-h>0$.

We know $z=1$ and $\eta = |B|/|A|^2 $ on the boundary of the disc of convergence in (\ref{eq:13}). 
Putting  (\ref{eq:5}) and  (\ref{eq:6}) at $\eta^2 \left|\overline{y}_2(z) \right|$ in (\ref{eq:14}) where $z=1$ and $\eta = |B|/|A|^2 $, we arrive at the following inequality such as
\begin{eqnarray}
\eta^2 \left|\overline{y}_2(z) \right|  
&>& \eta ^2 \sum_{i_0=0}^{\infty } \frac{(1-\epsilon )^{i_0}}{i_0+1+h_0+N} \sum_{i_2=i_0}^{\infty } \frac{(1-\epsilon )^{i_2 - i_0}}{i_2+3+h_0+N} \sum_{i_4=i_2}^{\infty } \frac{\left( 4+N-h +i_2\right)_{i_4 -i_2}}{\left( 4+N  +i_2\right)_{i_4 -i_2}}    \nonumber\\
&>& \eta ^2 \sum_{i_0=m}^{\infty } \frac{1}{i_0+1+h_0+N} \sum_{i_2=i_0}^{\infty } \frac{(1-\epsilon )^{i_2 }}{i_2+3+h_0+N} \sum_{i_4=i_2}^{\infty } \frac{\left( 4+N-h +i_2\right)_{i_4 }}{\left( 4+N  +i_2\right)_{i_4  }}
 \label{eq:19}
\end{eqnarray}
Put $r=1$ in (\ref{eq:16}) and take the new (\ref{eq:16}) into (\ref{eq:19})
\begin{eqnarray} 
\eta^2 \left|\overline{y}_2(z) \right|  
&>& \frac{\eta^2}{2} \sum_{i_0=m}^{\infty } \frac{1}{i_0+1+h_0+N}  \sum_{i_2=i_0}^{\infty }\frac{\Gamma \left( 4+N-h +i_2\right)}{\Gamma \left( 4+N +i_2\right)} \frac{(1-\epsilon )^{i_2  }}{i_2+3+h_0+N} \sum_{i_4=i_2}^{\infty } \frac{1}{i_4^{h}} \nonumber\\
&=&  \frac{\eta^2}{2}  \sum_{i_0=m}^{\infty } \frac{1}{i_0+1+h_0+N}  \sum_{i_2=m}^{\infty }\frac{\Gamma \left( 4+N-h +i_2\right)}{\Gamma \left( 4+N +i_2\right)} \frac{(1-\epsilon )^{i_2  }}{i_2+3+h_0+N} \sum_{i_4=m}^{\infty } \frac{1}{i_4^{h}} -\cdots  \nonumber\\
&=&  \frac{(1-\epsilon )^m \eta^2}{2}  \frac{\Gamma \left( N+m+4\right)}{(N+m+3+h_0)\Gamma \left( N+m+4-h\right)}  \nonumber\\
&&\times \pFq[4]{3}{2}{1, N\!+\!m\!+\!4, N\!+\!m\!+\!3\!+\!h_0}{N\!+\!m\!+\!4\!-\!h, N\!+\!m\!+\!4\!+\!h_0}{1\!-\!\epsilon}  \sum_{j=m}^{\infty } \frac{1}{N\!+\!1\!+\!j\!+\!h_0}\sum_{k=m}^{\infty } \frac{1}{ k^{h}}  -\cdots \hspace{1cm}
 \label{eq:20}
\end{eqnarray} 
Putting  (\ref{eq:5}) and  (\ref{eq:6}) at $\eta^3 \left|\overline{y}_3(z) \right|$ in (\ref{eq:14}) where $z=1$ and $\eta = |B|/|A|^2 $,  the inequality ensues such as
\begin{eqnarray}
\eta^3 \left|\overline{y}_3(z) \right|  
&>& \eta^3  \sum_{i_0=0}^{\infty } \frac{(1-\epsilon )^{i_0}}{i_0\!+\!1\!+\!h_0\!+\!N} \sum_{i_2=i_0}^{\infty } \frac{(1-\epsilon )^{i_2 - i_0}}{i_2\!+\!3\!+\!h_0\!+\!N} \sum_{i_4=i_2}^{\infty } \frac{(1-\epsilon )^{i_4 - i_2}}{i_4\!+\!5\!+\!h_0\!+\!N} \sum_{i_6=i_4}^{\infty } \frac{\left( 6\!+\!N\!-\!h\!+\!i_4\right)_{i_6 -i_4}}{\left(  6\!+\!N \!+\!i_4\right)_{i_6 -i_4}}      \nonumber\\
&>& \eta^3  \sum_{i_0=m}^{\infty } \frac{1}{  i_0\!+\!1\!+\!h_0\!+\!N } \sum_{i_2=i_0}^{\infty } \frac{1}{  i_2\!+\!3\!+\!h_0\!+\!N } \sum_{i_4=i_2}^{\infty } \frac{(1-\epsilon )^{i_4}}{  i_4\!+\!5\!+\!h_0\!+\!N } \sum_{i_6=i_4}^{\infty } \frac{\left(   6\!+\!N\!-\!h\!+\!i_4 \right)_{i_6 }}{\left(   6\!+\!N \!+\!i_4\right)_{i_6 }}   \hspace{1cm}
 \label{eq:21}
\end{eqnarray}  
Put $r=2$ in (\ref{eq:16}) and take the new (\ref{eq:16}) into (\ref{eq:21})
\begin{eqnarray}
\eta^3 \left|\overline{y}_3(z) \right|   
&>& \frac{\eta^3}{2}  \sum_{i_0=m}^{\infty }  \frac{1}{  i_0\!+\!1\!+\!h_0\!+\!N } \sum_{i_2=i_0}^{\infty } \frac{1}{  i_2\!+\!3\!+\!h_0\!+\!N } \sum_{i_4=i_2}^{\infty }\frac{\Gamma \left( 6\!+\!N\! +\!i_4\right)}{\Gamma \left( 6\!+\!N\!-\!h\! +\!i_4\right)} \frac{(1-\epsilon )^{i_4  }}{i_4 \!+\!5\!+\!h_0\!+\!N} \sum_{i_6=i_4}^{\infty }  \frac{1}{i_6^{h}} \nonumber\\
&=& \frac{\eta^3}{2}  \sum_{i_0=m}^{\infty }  \frac{1}{  i_0\!+\!1\!+\!h_0\!+\!N } \sum_{i_2=m}^{\infty } \frac{1}{  i_2\!+\!3\!+\!h_0\!+\!N } \sum_{i_4=m}^{\infty }\frac{\Gamma \left( 6\!+\!N\! +\!i_4\right)}{\Gamma \left( 6\!+\!N\!-\!h\! +\!i_4\right)} \frac{(1-\epsilon )^{i_4  }}{i_4 \!+\!5\!+\!h_0\!+\!N} \sum_{i_6=m}^{\infty }  \frac{1}{i_6^{h}}  -\cdots \nonumber\\
&=& \frac{(1-\epsilon )^m \eta^3}{2}  \frac{\Gamma \left( N+m+6\right)}{(N+m+5+h_0)\Gamma \left( N+m+6-h\right)}  \nonumber\\
&&\times \pFq[4]{3}{2}{1, N\!+\!m\!+\!6, N\!+\!m\!+\!5\!+\!h_0}{N\!+\!m\!+\!6\!-\!h, N\!+\!m\!+\!6\!+\!h_0}{1\!-\!\epsilon} \prod_{l=0}^{1}\left( \sum_{j=m}^{\infty } \frac{1}{N\!+\!2l\!+\!1\!+\!j\!+\!h_0} \right) \sum_{k=m}^{\infty } \frac{1}{ k^{h}}  -\cdots \hspace{1cm}
 \label{eq:22}
\end{eqnarray}
Putting  (\ref{eq:5}) and  (\ref{eq:6}) at $\eta^4 \left|\overline{y}_4(z) \right|$ in (\ref{eq:14}) where $z=1$ and $\eta = |B|/|A|^2 $, the  inequality is followed as
\begin{eqnarray}
\eta^4 \left|\overline{y}_4(z) \right|  
&>& \eta^4  \sum_{i_0=0}^{\infty } \frac{(1-\epsilon )^{i_0}}{i_0\!+\!1\!+\!h_0\!+\!N} \sum_{i_2=i_0}^{\infty } \frac{(1-\epsilon )^{i_2 - i_0}}{i_2\!+\!3\!+\!h_0\!+\!N} \sum_{i_4=i_2}^{\infty } \frac{(1-\epsilon )^{i_4 - i_2}}{i_4\!+\!5\!+\!h_0\!+\!N}  \nonumber\\
&&\times \sum_{i_6=i_4}^{\infty }\frac{(1-\epsilon )^{i_6 - i_4}}{i_6\!+\!7\!+\!h_0\!+\!N} \sum_{i_8=i_6}^{\infty } \frac{\left( 8\!+\!N\!-\!h\!+\!i_6\right)_{i_8 -i_6}}{\left(  8\!+\!N \!+\!i_6\right)_{i_8 -i_6}}      \nonumber\\
&>& \eta^4  \sum_{i_0=m}^{\infty } \frac{1}{i_0\!+\!1\!+\!h_0\!+\!N} \sum_{i_2=i_0}^{\infty } \frac{1}{i_2\!+\!3\!+\!h_0\!+\!N} \sum_{i_4=i_2}^{\infty } \frac{1}{i_4\!+\!5\!+\!h_0\!+\!N}  \nonumber\\
&&\times \sum_{i_6=i_4}^{\infty }\frac{(1-\epsilon )^{i_6 }}{i_6\!+\!7\!+\!h_0\!+\!N} \sum_{i_8=i_6}^{\infty } \frac{\left( 8\!+\!N\!-\!h\!+\!i_6\right)_{i_8  }}{\left(  8\!+\!N \!+\!i_6\right)_{i_8 }}  
 \label{eq:23}
\end{eqnarray} 
Put $r=3$ in (\ref{eq:16}) and take the new (\ref{eq:16}) into (\ref{eq:23})
\begin{eqnarray}
\eta^4 \left|\overline{y}_4(z) \right|   
&>& \frac{\eta^4}{2}  \sum_{i_0=m}^{\infty }  \frac{1}{  i_0\!+\!1\!+\!h_0\!+\!N } \sum_{i_2=i_0}^{\infty } \frac{1}{  i_2\!+\!3\!+\!h_0\!+\!N } 
\sum_{i_4=i_2}^{\infty } \frac{1}{  i_4\!+\!5\!+\!h_0\!+\!N } \nonumber\\
&&\times \sum_{i_6=i_4}^{\infty }  \frac{\Gamma \left( 8\!+\!N\! +\!i_6\right)}{\Gamma \left( 8\!+\!N\!-\!h\! +\!i_6\right)} \frac{(1-\epsilon )^{i_6 }}{i_6 \!+\!7\!+\!h_0\!+\!N} \sum_{i_8=i_6}^{\infty }  \frac{1}{i_8^{h}} \nonumber\\
&=& \frac{\eta^4}{2}  \sum_{i_0=m}^{\infty }  \frac{1}{  i_0\!+\!1\!+\!h_0\!+\!N } \sum_{i_2=m}^{\infty } \frac{1}{  i_2\!+\!3\!+\!h_0\!+\!N } 
\sum_{i_4=m}^{\infty } \frac{1}{  i_4\!+\!5\!+\!h_0\!+\!N } \nonumber\\
&&\times \sum_{i_6=m}^{\infty }  \frac{\Gamma \left( 8\!+\!N\! +\!i_6\right)}{\Gamma \left( 8\!+\!N\!-\!h\! +\!i_6\right)} \frac{(1-\epsilon )^{i_6 }}{i_6 \!+\!7\!+\!h_0\!+\!N} \sum_{i_8=m}^{\infty }  \frac{1}{i_8^{h}} -\cdots \nonumber\\
&=& \frac{(1-\epsilon )^m \eta^4}{2}  \frac{\Gamma \left( N+m+8\right)}{(N+m+7+h_0)\Gamma \left( N+m+8-h\right)}  \nonumber\\
&&\times \pFq[4]{3}{2}{1, N\!+\!m\!+\!8, N\!+\!m\!+\!7\!+\!h_0}{N\!+\!m\!+\!8\!-\!h, N\!+\!m\!+\!8\!+\!h_0}{1\!-\!\epsilon} \prod_{l=0}^{2}\left( \sum_{j=m}^{\infty } \frac{1}{N\!+\!2l\!+\!1\!+\!j\!+\!h_0} \right) \sum_{k=m}^{\infty } \frac{1}{ k^{h}}  -\cdots \hspace{1cm}
 \label{qq:24}
\end{eqnarray} 
By mathematical induction, we repeat this process and  construct inequalities of every $\eta^{\tau } \left|\overline{y}_{\tau }(z) \right|$ terms where $\tau  \geq  5$.
Substitute (\ref{eq:20}), (\ref{eq:22}), (\ref{qq:24}) and including  inequalities of all $\eta^{\tau } \left|\overline{y}_{\tau }(z) \right|$ terms where $\tau \geq  5$ into (\ref{eq:13})
\begin{eqnarray}
\sum_{i=0}^{\infty }\left|\overline{c}_{i}\right| |x|^{i} &>&  \frac{(1-\epsilon )^m }{2}  \sum_{p=1}^{\infty }\frac{\Gamma \left( N\!+\!m\!+\!2p\!+\!2\right)\eta ^{p+1}}{(N\!+\!m\!+\!2p\!+\!1\!+\!h_0)\Gamma \left( N\!+\!m\!+\!2p\!+\!2\!-\!h\right)}  \nonumber\\
&&\times \pFq[4]{3}{2}{1, N\!+\!m\!+\!2p\!+\!2, N\!+\!m\!+\!2p\!+\!1\!+\!h_0}{N\!+\!m\!+\!2p\!+\!2\!-\!h, N\!+\!m\!+\!2p\!+\!2\!+\!h_0}{1\!-\!\epsilon} \prod_{l=0}^{p-1}\left( \sum_{j=m}^{\infty } \frac{1}{N\!+\!2l\!+\!1\!+\!j\!+\!h_0} \right) \sum_{k=m}^{\infty } \frac{1}{ k^{h}}  -\cdots \nonumber\\
&>&  \frac{(1-K)(1-\epsilon )^m }{2}  \sum_{p=1}^{\infty }\frac{\Gamma \left( N\!+\!m\!+\!2p\!+\!2\right)\eta ^{p+1}}{(N\!+\!m\!+\!2p\!+\!1\!+\!h_0)\Gamma \left( N\!+\!m\!+\!2p\!+\!2\!-\!h\right)}  \nonumber\\
&&\times \pFq[4]{3}{2}{1, N\!+\!m\!+\!2p\!+\!2, N\!+\!m\!+\!2p\!+\!1\!+\!h_0}{N\!+\!m\!+\!2p\!+\!2\!-\!h, N\!+\!m\!+\!2p\!+\!2\!+\!h_0}{1\!-\!\epsilon} \prod_{l=0}^{p-1}\left( \sum_{j=m}^{\infty } \frac{1}{N\!+\!2l\!+\!1\!+\!j\!+\!h_0} \right) \sum_{k=m}^{\infty } \frac{1}{ k^{h}} \nonumber\\
&>&  \frac{(1\!-\!K)(1\!-\!\epsilon )^m }{2}  \sum_{p=1}^{\infty }\frac{\Gamma \left( N\!+\!m\!+\!2p\!+\!2\right)\eta ^{p+1}}{(N\!+\!m\!+\!2p\!+\!1\!+\!h_0)\Gamma \left( N\!+\!m\!+\!2p\!+\!2\!-\!h\right)}  \nonumber\\
&&\times \prod_{l=0}^{p-1}\left( \sum_{j=m}^{\infty } \frac{1}{N\!+\!2l\!+\!1\!+\!j\!+\!h_0} \right) \sum_{k=m}^{\infty } \frac{1}{ k^{h}} 
 \label{eq:24}
\end{eqnarray}
where $0<K<1$ and   $ \eta =|B|/|A|^2$. We know $\sum_{k=m}^{\infty } \frac{1}{ k^{h}} <\infty $ but a harmonic series $\sum_{j=m}^{\infty } \frac{1}{N+2l+1+j+h_0}>\infty $.  So, $\sum_{i=0}^{\infty }\left|\overline{c}_{i}\right| |x|^{i} >\infty $.  Similarly, an inequality of $\sum_{i=0}^{\infty }\left| \hat{c}_{i}\right| |x|^{i}$ is constructed by replacing $N$ with $N+1$ in (\ref{eq:24}). So, $\sum_{i=0}^{\infty }\left| \hat{c}_{i}\right| |x|^{i} >\infty$.
Therefore, (\ref{eq:12}) is divergent if $ \Omega _{t-1}< \omega _{t-1} $. 

\textbf{2. The case of} $\mathbf{\Omega _{t-1}\geq \omega _{t-1}}$   

If $\Omega _{t-1}\geq \omega _{t-1}$,
\begin{equation}
 \left| \overline{A}_n \right| >  1-\frac{h_1}{n}>1-\epsilon 
\label{eq:25}
\end{equation} 
for some positive integer $h_1$ with given positive error bound $\epsilon $ where $n\geq N$. The rigorous proof of (\ref{eq:25}) is available in Ref.\cite{Choun2020}.
(\ref{eq:5}) and (\ref{eq:25}) are in the same form of inequality.  We know that (\ref{eq:16}) is also satisfied with $h_1$. So, inequality of $\sum_{i=0}^{\infty }\left|\overline{c}_{i}\right| |x|^{i} $ is same as (\ref{eq:24}) by replacing $h$ with $h_1$ where $N-h_1>0$. So, $\sum_{i=0}^{\infty }\left|\overline{c}_{i}\right| |x|^{i} >\infty $ Similarly, $\sum_{i=0}^{\infty }\left| \hat{c}_{i}\right| |x|^{i} >\infty$. Therefore,  (\ref{eq:12}) is divergent if $ \Omega _{t-1}\geq \omega _{t-1} $.
  
Eventually, we conclude that  a confluent Heun function does not converge on the boundary of the disc of convergence. \qed
\end{pf} 

\begin{thm} 
Let $y(x)= x^{\lambda }\left( d_0 +d_1 x+ d_2 x^2+ d_3 x^3+ \cdots \right)$, $d_i \ne 0$, be a power series  which consists of the 3-term recurrence relation. We assume that the polynomials in the numerator and denominator of $B_n$ have the same degree, and the degree of polynomials in the numerator of $A_n$ is  smaller than the degree of polynomials in the denominator of it:
\begin{equation}
d_{n+1}=A_n \;d_n +B_n \;d_{n-1} \hspace{1cm};n\geq 1 \label{eq:100}
\end{equation} 
\begin{equation}
\begin{cases} 
A_n =  \dfrac{C_{t-1}n^{t-1}+ C_{t-2}n^{t-2} + \cdots + C_0}{c_t n^t + c_{t-1}n^{t-1} + \cdots + c_0}  = A  \dfrac{ n^{t-1} + \Omega _{t-2}n^{t-2} + \cdots + \Omega _0}{ n^t + \omega _{t-1}n^{t-1} + \cdots + \omega _0} \\
B_n = \dfrac{G_t n^t + G_{t-1}n^{t-1} + \cdots + G_0}{g_t n^t + g_{t-1}n^{t-1} + \cdots + g_0}  = B \dfrac{ n^t + \Theta _{t-1}n^{t-1} + \cdots + \Theta _0}{ n^t + \theta _{t-1}n^{t-1} + \cdots + \theta _0} \\
d_1= \dfrac{C_0}{c_0}  \;d_0 = A \dfrac{\Omega _0}{\omega _0}  \;d_0   
\end{cases}\nonumber
\end{equation}
where neither $C_{t-1}$ nor $c_t$ is zero and neither $G_{t}$ nor $g_t$ is zero with $\Omega _j= C_j/C_{t-1}$, $\omega _j = c_j/c_t $, $\Theta _j = G_j/G_{t} $ and $\theta _j= g_j/g_t$. Also, we denote $A=C_{t-1}/c_t$ and $B= G_{t}/g_t$.
 
The domain of absolute convergence of $y(x)$ is given by \cite{Choun2013} 
\begin{equation}
\mathcal{D} :=  \left\{  x \in \mathbb{C} \Bigg|  \left| \lim_{n\rightarrow \infty} A_n x \right| + \left| \lim_{n\rightarrow \infty} B_n x^2 \right| =  \left| B x^2 \right| <1 \right\}
 \nonumber
\end{equation}
If $ \left| B x^2 \right| =1$, then the series cannot converges. \label{Thm.2} 
\end{thm} 
\begin{pf}
\textbf{1. The case of}  $\mathbf{ \Theta _{t-1} < \theta _{t-1}}$ 

We define
\begin{equation} 
\left| \overline{A}_n \right| =  \left|  \dfrac{ n^{t-1} + \Omega _{t-2}n^{t-2} + \cdots + \Omega _0}{ n^t + \omega _{t-1}n^{t-1} + \cdots + \omega _0} \right|  \;\;\mbox{and}\;\; \left| \overline{B}_n \right| =  \left|  \dfrac{ n^t + \Theta _{t-1}n^{t-1} + \cdots + \Theta _0}{ n^t + \theta _{t-1}n^{t-1} + \cdots + \theta _0}  \right|    
 \nonumber
\end{equation} 
If  $\theta _{t-1} - \Theta _{t-1}- h_2 <0$ with a positive integer $h_2$ and given positive error bound $\epsilon $,    
\begin{equation}
 \left| \overline{B}_n \right| > 1-\frac{h_2}{n} > 1-\epsilon 
\label{eq:101}
\end{equation}
where $n\geq N$. 

And we say
\begin{equation}
 \left| \overline{A}_n \right| > \frac{1}{n+h_3}
\label{eq:102}
\end{equation}
for a large $n$ and some positive integer $h_3$.

(\ref{eq:14}) can described differently such as
\begin{align}
\left|\overline{y}_0(z)\right| &= \sum_{i_0=0}^{\infty }  \prod _{i_1=0}^{i_0-1}\left| \overline{B}_{2i_1+1+N} \right| z^{i_0}\nonumber\\
\left|\overline{y}_1(z) \right| &=  \sum_{i_0=0}^{\infty } \left| \overline{A}_{2i_0+N}\right| \prod _{i_1=0}^{i_0-1}\left|\overline{B}_{2i_1+1+N} \right| \sum_{i_2=i_0}^{\infty } \prod _{i_3=i_0}^{i_2-1}\left|\overline{B}_{2i_3+2+N}\right| z^{i_2}\nonumber\\
\left|\overline{y}_{\tau}(z)\right|  &=   \sum_{i_0=0}^{\infty } \left|\overline{A}_{2i_0+N}\right| \prod _{i_1=0}^{i_0-1} \left| \overline{B}_{2i_1+1+N}\right|
 \prod _{k=1}^{\tau -1} \left( \sum_{i_{2k}= i_{2(k-1)}}^{\infty } \left|\overline{A}_{2i_{2k}+k+N}\right| \prod _{i_{2k+1}=i_{2(k-1)}}^{i_{2k}-1}\left|\overline{B}_{2i_{2k+1}+k+1+N}\right|\right)  \nonumber\\
  & \hspace{.3cm}\times \sum_{i_{2\tau } = i_{2(\tau -1)}}^{\infty }  \prod _{i_{2\tau +1}=i_{2(\tau -1)}}^{i_{2\tau }-1} \left|\overline{B}_{2i_{2\tau +1}+\tau +1+N}\right|  z^{i_{2\tau }}    \label{eq:103}
\end{align}  
and
\begin{equation}
\begin{cases}
\eta = |A| |x| \cr
z= |B| |x|^2 
\end{cases}\nonumber 
\end{equation}
(\ref{eq:103}) is done by letting $\left|A_{\overline{n}}\right|$ in the sequence $\left|\overline{c}_{n}\right|$ is the leading term in a series $\sum_{i=0}^{\infty }\left|\overline{c}_{i}\right| |x|^{i}$; we observe the term of sequence $\left|\overline{c}_{n}\right|$ which includes zero term of $\left|A_{\overline{n}}\right|'s$ for a sub-power series $\left|\overline{y}_0(z)\right|$, one term of $\left|A_{\overline{n}}\right|'s$ for the sub-power series $\left|\overline{y}_1(z)\right|$, two terms of $\left|A_{\overline{n}}\right|'s$ for a $\left|\overline{y}_2(z)\right|$, three terms of $\left|A_{\overline{n}}\right|'s$ for a $\left|\overline{y}_3(z)\right|$, etc.

In the same way that we prove Thm.\ref{Thm.1}, as we put (\ref{eq:101}) and (\ref{eq:102}) in $\eta ^{\tau}\left|\overline{y}_{\tau}(z)\right| $ where $z=1$ and $\eta = |A|/\sqrt{|B|} $ at (\ref{eq:103}), harmonic series start to appear. So, $\sum_{i=0}^{\infty }\left|\overline{c}_{i}\right| |x|^{i}, \sum_{i=0}^{\infty }\left| \hat{c}_{i}\right| |x|^{i} >\infty $. Therefore,  (\ref{eq:12}) is divergent if $ \Theta _{t-1} < \theta _{t-1} $.

\textbf{2. The case of}  $\mathbf{ \Theta _{t-1} \geq \theta _{t-1}}$

If $  \Theta _{t-1} \geq \theta _{t-1} $,
 \begin{equation}
 \left| \overline{B}_n \right| > 1-\frac{h_4}{n} > 1-\epsilon 
\label{eq:104}
\end{equation}
for some positive integer $h_4$ with given positive error bound $\epsilon $ where $n\geq N$. 

By substituting (\ref{eq:104}) and (\ref{eq:102}) into $\eta ^{\tau}\left|\overline{y}_{\tau}(z)\right| $ where $z=1$ and $\eta = |A|/\sqrt{|B|} $ at (\ref{eq:103}), harmonic series are also shown. Then, $\sum_{i=0}^{\infty }\left|\overline{c}_{i}\right| |x|^{i}, \sum_{i=0}^{\infty }\left| \hat{c}_{i}\right| |x|^{i} >\infty $. Again,  (\ref{eq:12}) is divergent if $ \Theta _{t-1} \geq \theta _{t-1} $. 
\qed
\end{pf}

\end{document}